\documentclass[reqno,12pt]{amsart}

\title
[Lagrangian mean curvature flows]
{Generalized Lagrangian mean curvature flows in symplectic manifolds}
\author[Knut Smoczyk]{\sc Knut Smoczyk$^\ast$}
\address{
$^\ast$Leibniz Universit\"at Hannover,
{F}akult\"at f\"ur Mathematik und Physik,
Institut f\"ur Mathematik,
Welfengarten 1,
30167 Hannover,
Germany}
\email{smoczyk@math.uni-hannover.de}

\author[Mu-Tao Wang]{\sc Mu-Tao Wang$^{\ast\ast}$}
\address{
$^{\ast\ast}$Columbia University,
Department of Mathematics,
2990 Broadway,
New York, NY 10027}
\email{mtwang@math.columbia.edu}
\thanks{The first author was supported by the DFG (German Research Foundation). The second author was partially supported by National Science
Foundation Grant DMS 0605115 and 0904281.}%

\date{\today}

\usepackage{amscd,amsfonts,mathrsfs,amsthm,enumerate}
\usepackage{amssymb, amsmath}          
\usepackage[shortalphabetic]{amsrefs}    

\def\gmean  {{\overrightarrow{\widehat{\operatorname{H}}}}}
\def\contraction{\lrcorner}

\def\hn {{\widehat\nabla}}
\def\hT {{\widehat T}}
\def\hC {{\widehat\Gamma}}
\def\hR {{\widehat R}}
\def\hA {{\widehat A}}
\def\hh {{\widehat h}}
\def\hH {{\widehat H}}
\def\hr {{\widehat r}}
\def\hs {{\widehat s}}
\def\hrho{{\widehat\rho}}

\def\dd       #1#2#3{{#1}_{#2#3}}
\def\ddd      #1#2#3#4{{#1}_{#2#3#4}}

\def\uu       #1#2#3{{#1}^{#2#3}}

\def\du       #1#2#3{{#1}_{#2}^{\phantom{#2}{#3}}}
\def\udu      #1#2#3#4{{#1}^{{#2}\phantom{#3}{#4}}_{\phantom{#2}{#3}}}

\def\ddu      #1#2#3#4{{#1}_{#2#3}^{\phantom{#2#3}{#4}}}
\def\udd      #1#2#3#4{{#1}^{#2}_{\phantom{#2}{#3#4}}}

\newtheorem{thm}{Theorem}

\newtheorem{pro}{Proposition}
\newtheorem{lem}{Lemma}
\newtheorem{dfn}{Definition}
\newtheorem{exa}{Example}
\newtheorem{rem}{Remark}

\begin{document}
\begin{abstract}
An almost K\"ahler structure on a symplectic manifold $(N, \omega)$ consists of a Riemannian metric $g$ and an almost complex structure $J$ such that
the symplectic form $\omega$ satisfies $\omega(\cdot, \cdot)=g(J(\cdot), \cdot)$. Any symplectic manifold admits an almost K\"ahler structure and we refer to $(N, \omega, g, J)$ as an almost K\"ahler manifold. In this article, we propose a natural evolution equation to investigate the deformation of Lagrangian submanifolds
in almost K\"ahler manifolds. A metric and complex connection $\hn$ on $TN$ defines a
generalized mean curvature vector field  along any Lagrangian submanifold $M$ of $N$.
We study the evolution of $M$ along this vector field, which turns out to be a Lagrangian deformation,
as long as the connection $\hn$ satisfies an Einstein condition. This can be viewed as a generalization
of the classical Lagrangian mean curvature flow in K\"ahler-Einstein manifolds where the connection $\hn$ is the Levi-Civita connection of $g$.  Our result applies to the important case of Lagrangian submanifolds
in a cotangent bundle equipped with the canonical almost K\"ahler structure and to other generalization of Lagrangian mean curvature flows, such as the flow considered by
Behrndt \cite{b} in K\"ahler manifolds that are almost Einstein.
\end{abstract}
\maketitle

\section{Introduction}

Special Lagrangian submanifolds \cite{hl} and Lagrangian mean curvature flows
\cite{ty} attract much attentions due to their relations to the SYZ
conjecture \cite{syz} on mirror symmetry between Calabi-Yau manifolds. A Calabi-Yau,
or in general a K\"ahler-Einstein manifold, is a great place to study the
mean curvature flow as this process provides a Lagrangian deformation \cite{sm1}.
This important property  no  longer holds if the ambient space is a general symplectic manifold. However, there are important conjectures (see for example \cite{fss}) concerning the Lagrangian isotopy problem in general symplectic manifolds such as cotangent bundles which do not carry K\"ahler-Einstein structures.

In this article, we aim at defining a generalized Lagrangian mean curvature flow in general almost K\"ahler manifolds $N$. We consider generalized mean curvature vector fields $\gmean$ (see Definition \ref{gmcv} in \S 4) along Lagrangian
or more generally almost Lagrangian submanifolds, i.e. submanifolds $M$ for which $J(TM)\cap TM=\{0\}$.
The definition of the generalized mean curvature vector $\gmean$ relies on a choice of a complex and metric
connection $\hn$ on $TN$ that could carry non-trivial torsion $\hT$. We then say that a smooth family of
almost Lagrangian immersions
\[F:M \times [0, T)\rightarrow N\] satisfies the generalized mean curvature flow, if
\begin{equation}\label{gmcf}
\frac{\partial F}{\partial t}(p, t)=\gmean(p, t)\,,\quad\text{and}\quad F(M, 0)=M_0
\end{equation}
where $\gmean(p,t)$ is the generalized mean curvature vector of the almost Lagrangian submanifold
$M_t=F(M, t)$ at $F(p,t)$. This flow is uniquely defined up to tangential diffeomorphisms of $M$.

Recall that the Riemannian metric $g$ on $N$ also defines the
classical mean curvature vector  $\overrightarrow{\operatorname{H}}$ on $M$ through the first variation of volume.
$\gmean$ differs from $\overrightarrow{\operatorname{H}}$
by some lower order terms
involving the torsion of $\hn$. Since the mean curvature flow is a nonlinear
parabolic system which is non-degenerate after gauge fixing, the short time existence of the
generalized
mean curvature flow in the class of almost Lagrangian submanifolds can thus be established.

\begin{thm}\label{short_time}

Suppose $(N, \omega, g, J)$ is an almost K\"ahler manifold and $\hn$ is a complex and metric connection on $TN$.
For any initial smooth compact almost Lagrangian submanifold $M_0$, there
exists a maximal
time $T\in(0,\infty]$ so that the generalized mean curvature flow (\ref{gmcf}) exists smoothly
on $[0, T)$ in the class of almost Lagrangian submanifolds.
\end{thm}

We will show in Lemma \ref{mcv_relation} that the generalized mean curvature vector $\gmean$ is related to
a $1$-form $\hH$ on $M$ that can be seen as a generalization of the classical
mean curvature form (or Maslov form) of Lagrangian submanifolds of K\"ahler manifolds.

It turns out that $\hH$ is closed on Lagrangian submanifolds and the flow preserves the Lagrangian condition, if in addition the {\sl Ricci form} of $\widehat{\nabla}$,
$$\hrho(V,W):=\frac{1}{2}\operatorname{trace}(\hR(V,W)\circ J),$$ where $\hR$ is the curvature operator of $\hn$, satisfies the following Einstein condition:

\begin{dfn}
A metric and complex connection $\hn$ on an almost K\"ahler manifold $(N,\omega, g,J)$
is called {\sl Einstein}, if the Ricci form of $\hn$ satisfies
$$\hrho=f\omega$$
for some smooth function $f$ on $N$.
\end{dfn}

\begin{rem}
In general Einstein connections, if they exist, are not unique. For example the canonical
connection $\hn$
on the cotangent bundle $T^*M$ of a Riemannian
manifold $M$ equipped with the metric of type $I+III$ as defined in \cite{yi} is Einstein
(even
Ricci flat, i.e. the Ricci form vanishes). Any connection
$\tilde\nabla:=\hn+c\lambda\otimes J$,
where $\lambda$ is the Liouville form of $T^*M$ and $c$ some constant is also Einstein
with $f=-nc$.
\end{rem}

\begin{thm}\label{preserved}
Suppose $(N, \omega ,g,J)$ is an almost K\"ahler manifold and $\hn$ is a complex and metric
connection that satisfies the Einstein condition.
Suppose $M_0$ is a closed Lagrangian submanifold of $N$. Then the generalized mean curvature flow (\ref{gmcf}) with respect to $\hn$ preserves the Lagrangian
condition.
\end{thm}

This gives a new and large class of symplectic manifolds where Lagrangian mean curvature flows
can be defined. Our class includes the classical Lagrangian mean curvature flow in K\"ahler-Einstein
manifolds, the modified Lagrangian mean curvature flows considered by Behrndt \cite{b} for
K\"ahler manifolds that are almost Einstein and the important class of Lagrangian submanifolds
 in the cotangent bundle of any given Riemannian manifold.

Our strategy for proving that the Lagrangian condition is preserved is similar to that in \cite{sm1}.
We will first consider the flow in the larger (and open) class of almost Lagrangian submanifolds in
$N$. After the short time existence is established in this class, we apply the maximum principle to prove that the Lagrangian
condition is preserved.

\begin{rem}
Theorems \ref{short_time} and \ref{preserved} show that there exists a generalized Lagrangian mean curvature
flow on an almost K\"ahler manifold with an Einstein connection. It is an
interesting question to identify which almost K\"ahler manifold has this property. In \S 5, we give a list of examples that include all currently known cases to our knowledge.
\end{rem}

The article is organized as follows. In \S \ref{sec ak}, the geometry of almost K\"ahler manifolds
and the space of complex and metric connections are reviewed.
In \S \ref{sec lag} we consider the geometry of Lagrangian and almost Lagrangian
submanifolds of almost K\"ahler manifolds and define generalized second fundamental forms
and mean curvature forms with respect to a metric and complex connection $\widehat{\nabla}$. In particular we
show that the generalized mean curvature form with respect to an Einstein connection is closed on a Lagrangian submanifold.
In \S \ref{sec gmcf} we define the generalized mean curvature flow and prove our main Theorems
\ref{short_time} and \ref{preserved}. In \S 5, we present known examples of almost K\"ahler manifolds with Einstein connections.

{\sl Acknowledgments.}
The authors are grateful to Professor Mao-Pei Tsui, Professor Conan Leung and Dr. Lars Sch\"afer
for useful discussions.


\section{Almost K\"ahler manifolds}\label{sec ak}
In this section we recall some properties of almost K\"ahler manifolds.
\subsection{Almost K\"ahler manifolds and Einstein connections}~

Suppose $(N,\omega, g,J)$ is an almost K\"ahler manifold. This means that $(N,g)$ is a Riemannian manifold,
$J\in\operatorname{End}(TN)$ is an almost complex structure (i.e. an endomorphisms with
$J^2=-\operatorname{Id}$) and in
addition the characteristic $2$-form $\omega(\cdot, \cdot) :=g(J\cdot,\cdot)$ is symplectic, i.e.
$\omega$ is skew and $d\omega=0$. $J$ needs not be integrable and thus $(N,J)$ in general is
not a complex manifold. An almost K\"ahler manifold is K\"ahler if and only if $J$ is parallel
with respect to the Levi-Civita connection $\nabla$ of $g$.

We start with some general properties of almost K\"ahler manifolds. For an almost K\"ahler manifold $(N,\omega, g,J)$, let us define the class $\mathscr{C}$  of metric and
complex connections
$$\mathscr{C}:=\{\hn:\hn g=0,\,\hn J=0\}\,.$$
It is well known that this set is non-empty.
There exists a ``canonical" connection $\hn$ in $\mathscr{C}$
in the following sense: If $\nabla$ denotes the Levi-Civita connection of $g$, then $\hn$ defined by
$$\hn_XY:=\nabla_XY-\frac{1}{2}J(\nabla_X J)Y=\nabla_XY+\frac{1}{2}(\nabla_XJ)(JY)$$
lies in $\mathscr{C}$.

The following lemma can be easily proved:
\begin{lem}\label{lem connection}
If $\hn\in\mathscr{C}$ and $\sigma\in\Omega^1(N)$ is an arbitrary
$1$-form, then $\widetilde\nabla:=\hn+\sigma\otimes J\in\mathscr{C}$.
The curvature tensors $\widetilde R$ and $\hR$ of $\widetilde\nabla$ resp. $\hn$
are related by
$$\widetilde R=\hR+d\sigma\otimes J$$
and the torsions $\widetilde T$ and $\hT$ are related by
$$\widetilde T=\hT+\sigma\wedge J\,.$$
\end{lem}

Throughout this paper $\hn$ denotes a metric and complex connection. By the compatibility of $\omega$ and $J$,
we see that $\hn$ is also symplectic, i.e. $\hn\omega=0$. Since $\hn$ is metric, the curvature
tensor $\hR$ of $\hn$ satisfies
\begin{eqnarray}
\langle\hR(X,Y)V,W\rangle=-\langle \hR(Y,X)V,W\rangle=-\langle\hR(X,Y)W,V\rangle\,.\label{ak6}
\end{eqnarray}
Moreover, the fact that $\hn$ is complex implies
\begin{eqnarray}
\langle\hR(X,Y)JV,W\rangle=\langle J\hR(X,Y)V,W\rangle=\langle\hR(X,Y)JW,V\rangle\,.\label{ak7}
\end{eqnarray}
As $\langle\hR(X,Y)JV,W\rangle$ is symmetric in $V$ and $W$ we can take the trace
over these two arguments with respect to the Riemannian metric $g$ on $N$. This gives the
Ricci form of $\widehat{\nabla}$
$$\hrho(X,Y):=\frac{1}{2}\operatorname{trace}(\hR(X,Y)\circ J)\,.$$
The first Bianchi identity for $\hR$ is
\begin{eqnarray}
&\hR(X,Y)Z+\hR(Y,Z)X+\hR(Z,X)Y\nonumber\\
&=\hT(\hT(X,Y),Z)+\hT(\hT(Y,Z),X)+\hT(\hT(Z,X),Y)\nonumber\\
&+(\hn_X\hT)(Y,Z)+(\hn_Y\hT)(Z,X)+(\hn_Z\hT)(X,Y)\,,\label{ak5}
\end{eqnarray}
where
$$\hT(X,Y)=\hn_XY-\hn_YX-[X,Y]$$
is the torsion tensor of $\hn$.

Due to this weaker Bianchi identity, the Ricci form $\hrho$
in general does not satisfy the equation $\hrho(V,JW)=\hR ic(V,W)$, where
$\hR ic$ denotes the usual Ricci curvature
$\hR ic(V,W)=\operatorname{trace}\langle\hR(\cdot,V)\cdot,W\rangle$.

\section{The geometry of Lagrangian submanifolds of almost K\"ahler manifolds}\label{sec lag}

Suppose $\hn\in\mathscr{C}$ is an arbitrary metric and complex
connection on the tangent bundle $TN$ of an almost K\"ahler manifold $(N, \omega, g, J)$.

The torsion tensor $\hT^\alpha_{\beta\gamma}$ of $\hn$ is locally given by
$$\hT^\alpha_{\beta\gamma}:=\hC^\alpha_{\beta\gamma}-\hC^\alpha_{\gamma\beta}\,,$$
where $\hC^\alpha_{\beta\gamma}$ are the Christoffel symbols of $\hn$. Locally for a
vector field
$V=V^\alpha\partial/\partial y^\alpha\in\Gamma(TN)$ we have
$$\hn_\alpha V^\beta=V^\beta_{,\alpha}+\hC^\alpha_{\beta\gamma}V^\gamma\,.$$


Suppose now that $F:M\to N$ is a smooth immersion. The differential
$$dF=\frac{\partial F^\alpha}{\partial x^i}\frac{\partial}{\partial y ^\alpha}\otimes dx^i$$
is a smooth section in the bundle $F^{-1}TN\otimes T^*M$. Let $\nabla$ denote the
Levi-Civita connection
on $TM$ with respect to  the induced Riemannian metric $\dd gij=\dd g\alpha\beta F^\alpha_iF^\beta_j\,,$
where $F^\alpha_i:=\partial F^\alpha/\partial x^i$. A connection $\hn$ on $TN$
induces a connection $\hn^{E}$ on the pull-back bundle $E:=F^{-1}TN$ along $M$ by
$$\hn^{E}_X\sigma:=\hn_{dF(X)}\sigma\,,$$
where $X\in TM$ and $\sigma\in\Gamma(E)$.
A product connection $\hn^E\otimes\nabla$ on the bundle $E\otimes TM$ is then defined by
$$(\hn^E\otimes\nabla)_X(\sigma\otimes Y):=\hn^E_X\sigma\otimes Y+\sigma\otimes\nabla_XY\,.$$
Similarly we obtain product connections on all bundles of the form
$$E\otimes\underbrace{TM\otimes\dots\otimes TM}_{p-times}\otimes
\underbrace{T^*M\otimes\dots\otimes T^*M}_{q-times}\,.$$

In some abuse of notation, let us denote all these connections on bundles containing $E$
as a factor by $\hn^E$. Since on $E$ we also have the Levi-Civita connection induced by
$g$, there are always two different connections $\nabla$, $\hn^E$ on bundles containig $E$.
We shall use both.

For example, for a section
$$V^\alpha_i\frac{\partial}{\partial y^\alpha}\otimes dx^i\in\Gamma(F^{-1}TN\otimes T^*M)$$
we get
$$\hn^E_iV^\alpha_j=V^\alpha_{j,i}
-\Gamma^k_{ij}V^\alpha_k+\hC^\alpha_{\beta\gamma}F^\beta_iV^\gamma_j$$
and likewise
$$\nabla_iV^\alpha_j=V^\alpha_{j,i}
-\Gamma^k_{ij}V^\alpha_k+\Gamma^\alpha_{\beta\gamma}F^\beta_iV^\gamma_j\,,$$
where $\hC^\alpha_{\beta\gamma}$ resp. $\Gamma^\alpha_{\beta\gamma}$ are the Christoffel symbols
of $\hn$ resp. $\nabla$ on $TN$.

We derive
\begin{eqnarray}
\hn^E_i\hn^E_jV^\alpha_k-\hn^E_j\hn^E_iV^\alpha_k
&=&-R^m_{kij}V_m^\alpha
+\hR^\alpha_{\delta\beta\epsilon}F^\beta_i F^\epsilon_jV^\delta_k\,,\label{cod1}
\end{eqnarray}
where $R^m_{kij}$ denotes the Riemann curvature of the Levi-Civita connection on $TM$.

Let us define two second fundamental tensors
$$\hA^\alpha_{ij}:=\hn^E_iF^\alpha_j\,,\quad A^\alpha_{ij}:=\nabla_iF^\alpha_j\,.$$
We obtain
\begin{equation}\label{sec_sym}
\hA^\alpha_{ij}-\hA^\alpha_{ji}=(\hC^\alpha_{\beta\gamma}
-\hC^\alpha_{\gamma\beta})F^\beta_iF^\gamma_j=\hT^\alpha_{\beta\gamma}F^\beta_iF^\gamma_j
\end{equation}
and
$$A^\alpha_{ij}=\hA^\alpha_{ij}+(\Gamma^\alpha_{\beta\gamma}
-\hC^\alpha_{\beta\gamma})F^\beta_i F^\gamma_j\,.$$
Applying (\ref{cod1}) to the section $F^\alpha_k$ yields the Codazzi equation
\begin{eqnarray}
\hn^E_i\hA^\alpha_{jk}-\hn^E_j\hA^\alpha_{ik}
&=&-R^m_{kij}F_m^\alpha
+\hR^\alpha_{\delta\beta\epsilon}F^\beta_i F^\epsilon_jF^\delta_k\,.\label{cod2}
\end{eqnarray}
Note, that (\ref{cod2}) does not contain any torsion terms $\hT$ though the connection
$\hn$ on $TN$ in general has torsion. This is because the part of $\hn^E$ that acts on the
tangent bundle is given by the Levi-Civita connection. All information on $\hn$ contained
in the Codazzi equation is then encoded in the curvature term
$\hR^\alpha_{\delta\beta\epsilon}F^\beta_i F^\epsilon_jF^\delta_k$.

Suppose now that $F:M\to N$ is almost Lagrangian, i.e. $\dim M=n=\frac{1}{2}\dim N$
and the map
$$\phi:TM\to T^\perp M\,,\quad \phi V:=(JV)^\perp$$
provides an isomorphism between the tangent and the normal bundle $T^\perp M$ of $M$.
This holds if and only if the symmetric tensor
$$\eta(V,W):=\langle\phi V,\phi W\rangle$$
is invertible everywhere.

In local coordinates $(x^i)_{i=1,\dots,n}$ on $M$, $F^*\omega$ can be written in the form
$F^*\omega=\dd\omega ijdx^i\otimes dx^j$ with
$$\dd\omega ij:=\dd\omega\alpha\beta F^\alpha_i F^\beta_j\,.$$
For $\eta=\dd\eta ijdx^i\otimes dx^j$ we obtain
$$\dd\eta ij=\dd gij-\du\omega im\dd \omega jm\,,$$
where here and in the following indices will be raised and lowered by contraction with the metric
tensors $\uu gij$ and $\dd gij$ and the Einstein convention always applies. An exception will be
the inverse of $\dd\eta ij$ which we denote
by $\uu\eta ij$ (note that this differs from
$\uu gik\uu gjl\dd\eta kl$). From
$\uu gij=\uu gim\du\delta mj$ and $\du\delta mj=\dd\eta mk\uu\eta kj$ we observe
\begin{equation}\label{inverse}
\uu\eta ij=\uu gij+\uu\omega il\dd\omega sl\uu\eta sj\,.
\end{equation}

\begin{dfn}
Denote $$\hh_{kij}:=\langle\phi F_k,\hA_{ij}\rangle,$$
the generalized mean curvature form of $M$ is defined to be $\hH=\hH_idx^i$ where $$\hH_i:=\uu gkj \hh_{kij} =g^{kj}\langle (JF_k)^\perp, \widehat{\nabla}_i F_j\rangle $$
\end{dfn}

We also define two auxiliary tensors
$$\hr_{kij}:=\omega(F_k,\hA_{ij})\,, \text{and} \quad\hs_{kij}:=\langle F_k,\hA_{ij}\rangle.$$ The following relations can be easily verified:
\begin{equation}\label{H_i}
\hh_{kij}=\hr_{kij}-\du\omega km\hs_{mij}\,,\quad\hH_i
=\ddu\hr kik+\uu\omega mk\ddd\hs mik\,.
\end{equation}

In the rest of the section, we compute $d\hH$. To this end let us first compute $\nabla_l\hr_{kij}$. Here $\hr$
is considered as a section in $\Gamma(T^*M\otimes T^*M\otimes T^*M)$ and $\nabla$ denotes
the Levi-Civita connection on $M$.
From $\hr_{kij}=\dd\omega\alpha\beta F^\alpha_k\hA^\beta_{ij}$ we obtain by Leibniz' rule
for connections
\begin{eqnarray}
\nabla_l\hr_{kij}&=&\hn_\gamma\dd\omega\alpha\beta F^\gamma_lF^\alpha_k\hA^\beta_{ij}
+\dd\omega\alpha\beta\hA^\alpha_{lk}\hA^\beta_{ij}
+\dd\omega\alpha\beta F^\alpha_k\hn^E_l\hA^\beta_{ij}\nonumber\\
&\overset{(\hn\omega=0)}{=}&\dd\omega\alpha\beta\hA^\alpha_{lk}\hA^\beta_{ij}
+\dd\omega\alpha\beta F^\alpha_k\hn^E_l\hA^\beta_{ij}\,.\nonumber
\end{eqnarray}
Interchanging $l$ and $i$ and subtracting yields
\begin{eqnarray}
&&~\hspace{-70pt}\nabla_l\ddu\hr kik-\nabla_i\ddu\hr klk\nonumber\\
&=&\uu gkj\dd\omega\alpha\beta(\hA^\alpha_{lk}\hA^\beta_{ij}-\hA^\alpha_{ik}\hA^\beta_{lj})
+\uu gkj\dd\omega\alpha\beta F^\alpha_k(\hn^E_l\hA^\beta_{ij}-\hn^E_i\hA^\beta_{lj})
\nonumber\\
&\overset{(\ref{cod2})}{=}&2\uu gkj\dd\omega\alpha\beta
\hA^\alpha_{lk}\hA^\beta_{ij}
+\uu gkj\dd\omega\alpha\beta F^\alpha_k(- R^m_{jli}F^\beta_m
+\hR^\beta_{\gamma\delta\epsilon}F^\gamma_jF^\delta_lF^\epsilon_i),\nonumber
\end{eqnarray} which can be rewritten as
\begin{equation}\label{hat_r}
\nabla_l\ddu\hr kik-\nabla_i\ddu\hr klk=2\uu gkj\omega(\hA_{lk},\hA_{ij})
+\du\omega mjR^m_{jli}
+\uu gkj\omega(F_k,\hR(F_l,F_i)F_j)\,.
\end{equation}
Let us treat the first term.
From
\begin{eqnarray}
\hA_{lk}
&=&\langle\hA_{lk},F^m\rangle F_m+\langle \hA_{lk},\phi F_m\rangle\uu\eta mn\phi F_n\nonumber\\
&=&\udd\hs mlkF_m+\uu\eta mn(\ddd\hr mlk-\du\omega mp\ddd\hs plk)\phi F_n\nonumber
\end{eqnarray}
and
\begin{eqnarray}
\omega(F_m,\phi F_n)=\dd\eta mn\,,\quad\omega(\phi F_p,\phi F_q)=-\dd\omega pq
+\du\omega pk\du\omega ql\dd\omega kl
\end{eqnarray}
we obtain
\begin{eqnarray}
2\uu gkj\omega(\hA_{lk},\hA_{ij})
&=&2\omega(\udd\hs mlkF_m+\uu\eta mn(\ddd\hr mlk-\du\omega mp\ddd\hs plk)\phi F_n,\nonumber\\
&&\udu\hs uikF_u+\uu\eta uv(\ddu\hr uik-\du\omega uq\ddu\hs qik)\phi F_v)\nonumber\\
&=&2\udd\hs mlk\uu\eta uv\ddu\hr uik\dd\eta mv-2\uu\eta mn\hr_{mlk}\udu\hs uik\dd\eta nu
+(C_1\contraction F^*\omega)_{li}\nonumber\\
&=&2\ddd\hs mlk\udu\hr mik-2\ddd\hs mik\udu\hr mlk+(C_1\contraction F^*\omega)_{li}\,,\nonumber
\end{eqnarray}
for some tensor $C_1$, where here and in the following $C\contraction F^*\omega$
denotes any tensor that is formed
by contracting an arbitrary tensor $C$ with $F^*\omega$.
Since $\hn$ is metric, the skew symmetry
\begin{equation}\label{sym s}
\hs_{kij}=-\hs_{jik}\,
\end{equation} holds.
Moreover, since $\hn$ is also symplectic applying Leibniz' rule with the
Levi-Civita connection $\nabla$ acting on $\omega$, we obtain
\begin{equation}\label{sym r}
\nabla_k\dd\omega ij=\omega(\hA_{ki},F_j)+\omega(F_i,\hA_{kj})=
\ddd\hr ikj-\ddd \hr jki\,.
\end{equation}
This implies
\begin{eqnarray}
2\ddd\hs mlk\udu\hr mik
&=&\ddd\hs mlk(\udu\hr mik-\udu \hr kim)\nonumber\\
&=&\ddd\hs mlk\nabla_i\uu\omega mk\nonumber\\
&=&\nabla_i(\ddd\hs mlk\uu\omega mk)-\nabla_i\ddd\hs mlk\uu\omega mk\,.\nonumber
\end{eqnarray}
Hence we conclude
\begin{eqnarray}\label{second_fund}
2\uu gkj\omega(\hA_{lk},\hA_{ij})
&=&\nabla_i(\ddd\hs mlk\uu\omega mk)-\nabla_l(\ddd\hs mik\uu\omega mk)
+(C_2\contraction F^*\omega)_{li}\,
\end{eqnarray}
where $C_2\contraction F^*\omega$ contains
the term $\nabla_l\ddd\hs mik\uu\omega mk-\nabla_i\ddd\hs mlk\uu\omega mk$.

With these preparations, we are ready to prove the main proposition of this section:
\begin{pro}
Let $(N,\omega, g,J)$ be an almost K\"ahler manifold and $\hn$ a metric and complex
connection on $TN$. Then the generalized mean curvature form $\hH$ for an
almost Lagrangian smooth immersion $F:M\to N$ satisfies
\begin{equation}\label{dh2}
d\hH=(C\contraction F^*\omega)-F^*\hrho
\end{equation}
for some smooth tensor field $C$. In particular, if $F$ is Lagrangian and
$\hn$ Einstein, then $\hH$ is closed.
\end{pro}

\begin{proof}
Plug equation (\ref{second_fund}) into equation (\ref{hat_r}) and recall equation (\ref{H_i}), we derive
\begin{eqnarray}
\nabla_l\hH_i-\nabla_i\hH_l
&=&\nabla_l\ddu\hr kik-\nabla_i\ddu\hr klk+\nabla_l(\ddd\hs mik\uu\omega mk)
-\nabla_i(\ddd\hs mlk\uu\omega mk)\nonumber\\
&=&(C_3\contraction F^*\omega)_{li}+\uu gkj\omega(F_k,\hR(F_l,F_i)F_j)\,.\label{dh1}
\end{eqnarray}

Since $\hn$ is metric and complex, we know that
$$\omega(\hR(V,W)X,Y)=\omega(\hR(V,W)Y,X)\,\text{for any}\,\, V,W,X,Y\in TN\,.$$
The Ricci form $\hrho$ is given by
$$\hrho(V,W):=\frac{1}{2}\sum_{\alpha=1}^{2n}\omega(\hR(V,W)e_\alpha,e_\alpha)\,,$$
where $e_\alpha$ is an arbitrary orthonormal basis of $TN$. In terms of $F_i$ and
$\phi F_i$, we can rewrite the trace in the form
$$\hrho(V,W)=\frac{1}{2}\uu gij\omega(\hR(V,W)F_i,F_j)+\frac{1}{2}\uu\eta ij\omega(\hR(V,W)\phi F_i,\phi F_j)\,.$$
On the other hand,
\begin{eqnarray}
\omega(\hR(V,W)\phi F_i,\phi F_j)
&=&\omega(\hR(V,W)(JF_i-\du\omega ilF_l),JF_j-\du\omega jkF_k)\nonumber\\
&=&\omega(J\hR(V,W)F_i,JF_j)+(C_4\contraction F^*\omega)(V,W)\nonumber\\
&=&\omega(\hR(V,W)F_i,F_j)+(C_4\contraction F^*\omega)(V,W)\nonumber
\end{eqnarray}
so that with (\ref{inverse}) we obtain
\begin{equation}\label{curvrel}
\hrho(V,W)=\uu gij\omega(\hR(V,W)F_i,F_j)+(C_5\contraction F^*\omega)(V,W)\,.
\end{equation}
Combining (\ref{dh1}) with (\ref{curvrel}) and taking into account that the Levi-Civita
connection $\nabla$ is torsion free
(so that $(d\hH)_{li}=\nabla_l\hH_i-\nabla_i\hH_l$), the proposition is proved.

\end{proof}


\section{The generalized mean curvature flow in almost K\"ahler manifolds}\label{sec gmcf}

\subsection{The generalized mean curvature vector}~
Suppose $F:M\to N$ is an almost Lagrangian submanifold of an almost K\"ahler manifold
$(N,\omega, g,J)$.
We recall that, with respect to the variation of volume defined by $g$, there is the classical mean curvature vector field $\overrightarrow{\operatorname{H}}$ defined on $M$:\[\overrightarrow{\operatorname{H}}=\eta^{kl} (A_{kl})^\perp=\eta^{ij} g^{kl}\langle A_{kl}, \phi F_i\rangle \phi F_j.\]
\begin{dfn}\label{gmcv}
 Let $\hn$ be a metric and complex connection.
The generalized mean curvature
vector of $M$ with respect to $\hn$ is defined to be
\begin{eqnarray}
\gmean
&:=&\overrightarrow{\operatorname{H}}+\uu \eta ij\uu gkl\left(\langle\hT(\phi F_i,F_k),F_l\rangle
+\langle\hT(F_i,F_k),\phi F_l\rangle\right)\phi F_j\,.\nonumber
\end{eqnarray}
\end{dfn}

In the following, we derive a relation between
the generalized mean curvature vector and the
generalized mean curvature form $\hH=\hH_idx^i$.

First we express the difference between $\hn$ and the Levi-Civita connection $\nabla$ of
$\langle\cdot, \cdot\rangle$ in terms of the torsion $\hT$ in the next Lemma.
\begin{lem}\label{connection_relation}
The two connections $\hn$ and $\nabla$ on $TN$ are related by
$$2\langle \hn_X Y-\nabla_X Y, Z\rangle
=\langle \hT(X, Y), Z\rangle+\langle \hT(Z, X), Y\rangle
+\langle \hT(Z, Y), X\rangle.$$
\end{lem}
\begin{proof}
Using the  fact that $\nabla$ is torsion free and compatible with the metric, we have
\begin{eqnarray}
&&X\langle Y, Z\rangle+Y\langle X, Z\rangle-Z\langle X, Y\rangle\nonumber\\
&&=\langle 2\nabla_X Y-[X, Y], Z\rangle+\langle Y, [X, Z]\rangle+\langle X, [Y, Z]\rangle.\nonumber
\end{eqnarray}
On the other hand, $\widehat{\nabla}$ is a metric connection, therefore
\begin{eqnarray}
&&X\langle Y, Z\rangle+Y\langle X, Z\rangle-Z\langle X, Y\rangle=
\langle 2\widehat{\nabla}_X Y-[X, Y]-\widehat{T}(X, Y), Z\rangle\nonumber\\
&&+\langle Y, [X, Z]+\widehat{T}(X, Z) \rangle+\langle X, [Y, Z]+\widehat{T}(Y, Z)\rangle\,.\nonumber
\end{eqnarray}
Subtracting the two identities yields the desired formula.
\end{proof}

\begin{lem}\label{mcv_relation}
The generalized mean curvature vector
$\gmean$ and the generalized mean curvature form $\hH$ with respect to the connection $\hn$ are related by
\begin{eqnarray}
\gmean=\uu\eta ij\left(\hH_i-\nabla^k\dd\omega ki-\du\omega im\ddu\hs mkk-\uu\omega mk\ddd\hs mki\right)\phi F_j\,.\label{eq gmeanrel}
\end{eqnarray}
In particular, if $M$ is Lagrangian, then
\begin{equation}\label{gmeanlag}
\gmean=\uu gij \hH_i J F_j\,.
\end{equation}
\end{lem}
\begin{proof}
From (\ref{H_i}), (\ref{sym s}) and (\ref{sym r}) we conclude
\begin{equation}\label{rel h1}
\ddd\hh lki=\ddd\hh ikl+\nabla_k\dd\omega li-\du\omega lm\ddd\hs mki+\du\omega im\ddd\hs mkl\,.
\end{equation}
Moreover, since $\hn$ has torsion we get
\begin{equation}\label{rel h2}
\ddd\hh lik=\ddd\hh lki+\langle\hT(F_i,F_k),\phi F_l\rangle\,.
\end{equation}
From these equations we then deduce
\begin{equation}\label{rel h3}
\ddu\hh ikk+\uu gkl\langle\hT(F_i,F_k),\phi F_l\rangle
=\hH_i-\nabla^k\dd\omega ki-\du\omega im\ddu\hs mkk-\uu\omega mk\ddd\hs mki\,.
\end{equation}
In addition, Lemma \ref{connection_relation} implies
\begin{eqnarray}
\ddu\hh ikk&=&\langle\du\hA kk,\phi F_i\rangle\nonumber\\
&=&\langle\overrightarrow H,\phi F_i\rangle+\uu gkl\langle\hT(\phi F_i,F_k),F_l\rangle\,.\nonumber
\end{eqnarray}
Combining this with (\ref{rel h3}) we finally get
\begin{eqnarray}
\langle\gmean,\phi F_i\rangle
&=&\langle \overrightarrow H,\phi F_i\rangle+\uu gkl\langle\hT(\phi F_i,F_k),F_l\rangle+\uu gkl\langle\hT(F_i,F_k),\phi F_l\rangle\nonumber\\
&=&\hH_i-\nabla^k\dd\omega ki-\du\omega im\ddu\hs mkk-\uu\omega mk\ddd\hs mki\,.\nonumber
\end{eqnarray}
This proves the lemma.
\end{proof}


\subsection{Proof of main theorems}~

We recall Theorem \ref{short_time}:
\setcounter{thm}{0}
\begin{thm}
Suppose $(N, \omega, g, J)$ is an almost K\"ahler manifold and $\hn$ is a complex and metric connection on $TN$.
For any initial smooth compact almost Lagrangian submanifold $M_0$, there
exists a maximal
time $T\in(0,\infty]$ so that the generalized mean curvature flow (\ref{gmcf}) exists smoothly
on $[0, T)$ in the class of almost Lagrangian submanifolds.
\end{thm}
\begin{proof}
It suffices to prove that the operator $E[F]:=\gmean[F]$ has no non-trivial degeneracies, i.e. we have to show that it is elliptic in the normal directions.
By Definition 3, the generalized mean curvature vector differs from the classical mean curvature vector only by terms of lower order so that the symbol of our operator is the same
as for the mean curvature flow and short-time existence follows.
\end{proof}

Next we recall Theorem \ref{preserved}:
\setcounter{thm}{1}
\begin{thm}
Suppose $(N,\omega,g,J)$ is an almost K\"ahler manifold and suppose $\hn$ is a complex and metric
connection that satisfies the Einstein condition.
Then the generalized mean curvature flow (\ref{gmcf}) w.r.t. $\hn$ preserves the Lagrangian
condition.
\end{thm}
\begin{proof}
From Cartan's formula we know
$$\frac{\partial}{\partial t}\,F^*\omega
=d\left(F^*\left(\frac{\partial F}{\partial t}\contraction\omega\right)\right)=d(F^*(\gmean\contraction\omega))\,.$$
We fix some time interval $[0,t_0]$, $0<t_0<T$.
From (\ref{dh2}), the Einstein property of $\hn$, and
(\ref{eq gmeanrel}) we obtain
\begin{eqnarray}
\frac{\partial}{\partial t}\,F^*\omega
&=&-d\hH+dd^\dagger (F^*\omega)+C_1\contraction F^*\omega+C_2\contraction\nabla( F^*\omega)\nonumber\\
&=&dd^\dagger (F^*\omega)+C_3\contraction F^*\omega+C_2\contraction\nabla(F^*\omega)\nonumber
\end{eqnarray}
for smooth tensor fields $C_1, C_2, C_3$, where
$d^\dagger(F^*\omega)$ is the $1$-form $\nabla^k\dd\omega kidx^i$.
Since $\omega$ (and $F^*\omega$) is closed we have
$$\Delta(F^*\omega)=dd^\dagger(F^*\omega)+C_4\contraction F^*\omega\,,$$
where $C_4$ depends on the Riemannian curvature of $M$. Combining the last identities we deduce
$$\frac{\partial}{\partial t}\,F^*\omega=\Delta(F^*\omega)+C_5\contraction F^*\omega+C_2\contraction \nabla(F^*\omega)$$
for smooth tensor fields $C_2, C_5$.

Together with
Cauchy-Schwarz' inequality, we thus obtain an estimate of the form
$$\frac{\partial}{\partial t}|F^*\omega|^2\le \Delta|F^*\omega|^2+c|F^*\omega|^2$$
for all $t\in[0,t_0]$ and some constant $c$ depending on $t_0$.

 Thus the growth rate of $\displaystyle|F^{*}\omega|^2$ on $[0,t_0]$ is at most exponential, i.e.
$$\sup_{p\in M}|F^*\omega|^2(p,t)\le \sup_{p\in M}|F^*\omega|^2(p,0)e^{ct}\,,\forall t\in[0,t_0].$$
However, $|F^*\omega|^2(p,0)$ is zero for all $p\in M$ as $M_0$ is Lagrangian.
Since $t_0$ is arbitrary, the theorem follows.
\end{proof}

The definition of  generalized Lagrangian mean curvature flows only
requires the existence of an ``almost Einstein connection" on $N$
in the sense that $\hrho-f\omega$ is $dd^c$ exact for some smooth
function $f$. In fact, the latter condition implies the existence
of an actual ``Einstein connection" on $N$, see Example 2 in the next section.

\section{Examples of almost K\"ahler manifolds with Einstein connections}

\begin{exa}

If $(N,\omega, g,J)$ is K\"ahler-Einstein, then we can choose $\hn$ to be the Levi-Civita connection
$\nabla$ of $g$ and $f=K/2n$, where $K$ is the scalar curvature of $N$ and $n$ is the complex dimension of $N$. In this case we recover
the classical Lagrangian mean
curvature flow.

\end{exa}
\begin{exa}
Let $(N,\omega, g,J)$ be one of the K\"ahler manifolds that are almost Einstein considered
in \cite{b}. In this case the Levi-Civita connection $\nabla$ is metric and complex and its
Ricci form
$\rho$ satisfies
$$\rho=\lambda\omega+ndd^c\psi$$
for some constant $\lambda$ and some smooth function $\psi$.  According to Lemma \ref{lem connection}
the connections
$$\hn:=\nabla+\sigma\otimes J$$
for any $\sigma\in\Omega^1(N)$
are also complex and metric and the curvature tensors $R$, $\hR$ of $\nabla$ resp.  $\hn$ are
related by $\hR=R+d\sigma\otimes J$. For the Ricci forms we get
$$\hrho=\rho-nd\sigma\,.$$
Therefore, if we choose $\sigma:=d^c\psi$, then the Ricci form $\hrho$ is conformal to $\omega$
(with $f=\lambda$). Since the torsion of $\hn$ is given by $\hT=\sigma\wedge J$, Proposition
\ref{mcv_relation} shows that the mean curvature vector
$\gmean$ w.r.t. the connection $\hn$ coincides with the mean curvature vector considered in
\cite{b} and we obtain the same flow. Since we do not need the integrability of $J$, we note that
the same trick works for almost K\"ahler manifolds
that are almost Einstein, so that the K\"ahler condition in \cite{b} is actually not needed.
\end{exa}

\begin{exa}
Let $N:=T^*M$ be the cotangent bundle of a Riemannian
manifold $(M,g)$.
The cotangent bundle carries a natural almost K\"ahler structure and a complex and metric
connection $\hn$ completely determined by the Levi-Civita connection
$\nabla$ of $(M,g)$
(see \cite{va} and \cite{yi}, \S IV.6). In this case the Ricci form of
the connection $\hn$
even vanishes (i.e. $f=0$) and our theorem applies. We will treat this example in great
detail in a forthcoming paper.

\end{exa}

\end{document}